\documentclass[12pt, reqno]{amsart}
\usepackage{amsmath, amsthm, amscd, amsfonts, amssymb, graphicx, color}
\usepackage{setspace}
\usepackage{mathrsfs}
\usepackage{multicol}
\usepackage[bookmarksnumbered, colorlinks, plainpages]{hyperref}
\hypersetup{colorlinks=true,linkcolor=red, anchorcolor=green, citecolor=cyan, urlcolor=red, filecolor=magenta, pdftoolbar=true}

\newtheorem{theorem}{Theorem}[section]
\newtheorem{lemma}[theorem]{Lemma}

\theoremstyle{definition}

\theoremstyle{remark}

\numberwithin{equation}{section}

\newcommand{\CC}{\mathbb {C}}

\begin{document}
\setcounter{page}{1}
\title[On the spectrum of Volterra-type  integral  operators on Fock--Sobolev spaces ] {On the spectrum of Volterra-type  integral operators on Fock--Sobolev spaces}
\author [Tesfa  Mengestie]{Tesfa  Mengestie }
\address{Department of Mathematical Sciences \\
Stord/Haugesund University College (HSH)\\
Klingenbergvegen 8, N-5414 Stord, Norway}
\email{Tesfa.Mengestie@hsh.no}
\thanks{The  author  is partially  supported by HSH grant 1244/H15.}
\subjclass[2010]{Primary 47B32, 30H20; Secondary 46E22, 46E20, 47B33}
 \keywords{Fock--Sobolev spaces, Weighted Fock spaces, Spectrum, Bounded, Compact, Volterra-type integral operators}
\begin{abstract}
We  determine the spectrum of the  Voltterra-type integral operators $V_g$  on the growth type Fock--Sobolev spaces $\mathcal{F}_{\psi_m}^\infty$. We also characterized the bounded and compact spectral properties of the operators in terms of function-theoretic properties of the inducing map $g$.     As a means to prove our main results, we first described  the spaces in terms of Littlewood--Paley type  formula which is interest of its own.
\end{abstract}

\maketitle

\section{Introduction} \label{1}
Let $m$  be  a nonnegative  integer and $0<p\leq \infty$. Then the Fock--Sobolev spaces  $\mathcal{F}_{(m,p)}$ consist of entire functions $f$ such that $f^{(m)}$, the $m$-th order derivative of $f$, belongs to the classical Fock spaces $\mathcal{F}_p;$
 which consist of all entire functions $f$ for which
\begin{align*}
\int_{\CC}
|f(z)|^p  e^{-\frac{p}{2}|z|^2} dA(z) <\infty  \ \ \text{and}\ \ \sup_{z\in \CC} |f(z)|e^{-\frac{1}{2}|z|^2}<\infty,
\end{align*}  respectively for  finite  and   infinite values of the exponents $p$, and  $dA$  here denotes the
usual Lebesgue area  measure on $\CC$.
 The Fock--Sobolev spaces were introduced in \cite{RCKZ} and studied further by several authors in different contexts   for example \cite{CCK, TM4, TM5}. Because of their Fourier  characterizations in \cite[Theorem A]{RCKZ}, the spaces can be simply considered as weighted Fock spaces induced by  the sequence of  weight functions $$\psi_m(z)= \frac{1}{2}|z|^2-m\log(1+ |z|). $$
 In  view of this,  $\mathcal{F}_{(m, p)}$ are just the weighted Fock spaces   $\mathcal{F}_{\psi_m}^p$ which consist of all entire functions $f$ for which \footnote{The  notation $U(z)\lesssim V(z)$ (or
equivalently $V(z)\gtrsim U(z)$) means that there is a constant
$C$ such that $U(z)\leq CV(z)$ holds for all $z$ in the set of a
question. We write $U(z)\simeq V(z)$ if both $U(z)\lesssim V(z)$
and $V(z)\lesssim U(z)$.}
\begin{align*}
%\label{normequal}
 \int_{\CC}
|f(z)|^p  e^{-p\psi_m(z) } dA(z)\simeq  \|f\|_{(m, p)}^p <\infty, \end{align*}
  for $0<p<\infty$ and  for $p=\infty$, the corresponding  estimate becomes
  \begin{align*}
 % \label{infinite}
  \|f\|_{(m,\infty)} \simeq \sup_{z\in \CC} |f(z)|e^{-\psi_m(z)} <\infty.
  \end{align*}
See  \cite{TM5} for further analysis on this.

The theory of integral operators constitutes a significant part of modern
  functional analysis, see for example  \cite{3CC,2KD,4DM,10HM} for some overviews.  A typical examples of these operators  includes the integral operators of Volterra. In this paper, we study  some spectral properties of  linear integral operators of Volterra-type.  More precisely,  for  a holomorphic function $g$,  we consider  the  Volterra-type
integral operator $V_g$ defined by
\begin{align*}
 V_gf(z)= \int_0^z f(w)g'(w) dw.
 \end{align*}
 These operators have been studied  for decades  in the settings of   various analytic function  spaces. We
have no intention to review a vast literature about it, but  for the state of the art, one can see \cite{Alman, Olivia1, Olivia,JPP, TM3, TM0, Si} and the related references
therein. Very recently,  some spectral properties of these operators were investigated in \cite{TM5} when they  act between the  Fock--Soblolev spaces $\mathcal{F}_{\psi_m}^p$ and $\mathcal{F}_{\psi_m}^q$  in which  both the exponents $p$ and $q$ are set to be finite. The central aim of this note is to continue that line of research and study the properties when at least one of the exponents is at infinity. Because of the monotonicity property of the spaces in the sense  $\mathcal{F}_{\psi_m}^p \subseteq \mathcal{F}_{\psi_m}^\infty$  for all exponent $0<p \leq \infty,$ one would expect that  a fairly weaker condition on the symbol $g$ in contrast to the finite exponent cases can give a bounded (compact) $V_g: \mathcal{F}_{\psi_m}^p \to \mathcal{F}_{\psi_m}^\infty$.  As can be seen in the  next section, it turns out that this is not the case, and the conditions are the same for both the finite and the infinite exponent cases  in dependent of the order $m$ and the size of the exponents as long the exponent in the target spaces is at least as big as in the domain space.  More precisely, for $0<p\leq \infty,$ it is proved that $V_g: \mathcal{F}_{\psi_m}^p \to \mathcal{F}_{\psi_m}^\infty$ is bounded if and only if $g$ is a complex polynomial of degree not bigger than two, while its compactness is characterized in terms of degree $g$ not being bigger than one. On the other hand, it is shown that  $V_g: \mathcal{F}_{\psi_m}^\infty\to \mathcal{F}_{\psi_m}^p$ is bounded (compact) if and only if  $g$ is of degree at most one.

Next we recall the notion of spectrum.  The spectrum $\sigma(T)$ of  a bounded operator $T$ on a Banach space consists of all $\lambda \in \CC$ for which $\lambda I- T$ is not invertible, where $I$ is the identity operator.  If $\lambda$ is an eigenvalue of $T,$ then the operator $\lambda I- T$ fails to be one to one and hence $\lambda I- T$ does not have inverse. The set of all such eigenvalues is referred to as the point spectrum of T and denoted by $\sigma_p(T)$:
$ \sigma_p(T)= \{\lambda \in \CC: Ker(\lambda I- T) \neq 0\}.$
 It  follows from this that $\sigma_p(T)\subseteq \sigma(T).$

In contrast to the boundedness, compactness,  and Schatten class membership  spectral properties of the integral operators $V_g$, there has not been  much studies on  their spectra. Some recent results in this connection can be read in \cite{Olivia2,TM5}. One of our  main  results describes the spectrum of $V_g$ acting on $\mathcal{F}_{\psi_m}^\infty$ in terms of a closed disk of center at the origin and radius involving the coefficient of the highest degree  term in  a polynomial expansion of $g$ as precisely stated   in our main result below.
\begin{theorem}\label{thm1}
Let $g$ be an entire function on $\CC$ and $0<p\leq  \infty$. Then
\begin{enumerate}
\item  $V_g: \mathcal{F}_{\psi_m}^p \to \mathcal{F}_{\psi_m}^\infty$ is
\begin{enumerate}
\item   bounded if and only if $g (z)= az^2+bz+c,\ \  a, b, c\in \CC$.% is a complex polynomial of degree not exceeding two.
\item  compact if and only if $g (z)= az+b, \ \  a, b \in \CC $.% is a complex polynomial of degree not exceeding one.
\end{enumerate}
\item if $0< p< \infty$, then  the following statements are equivalent:
\begin{enumerate}
\item $V_g: \mathcal{F}_{\psi_m}^\infty \to \mathcal{F}_{\psi_m}^p$ is   bounded;
\item $V_g: \mathcal{F}_{\psi_m}^\infty \to \mathcal{F}_{\psi_m}^p$ is   compact;
\item $g(z)= az+b, \ \ a, b\in \CC$ whenever $p>2$, and $g=$ constant otherwise.
\end{enumerate}
\item if $V_g:\mathcal{F}_{\psi_m}^\infty \to \mathcal{F}_{\psi_m}^\infty$  is   a bounded operator, i.e.  $g(z)= az^2+bz+c,\ \   a, b, c \in \CC$, then  we have
    \begin{align*}
\sigma(V_g)= \{\lambda \in \CC: |\lambda|\leq 2|a|\}= \{0\}\cup \overline{ \{\lambda \in \CC\setminus\{0\}:e^{g(z)/\lambda} \notin \mathcal{F}_{\psi_m}^\infty\} }.
\end{align*}
\end{enumerate}
\end{theorem}
As pointed earlier, when the operator $V_g$ acts between the Fock-Sobolev spaces $\mathcal{F}_{\psi_m}^p $ and $ \mathcal{F}_{\psi_m}^q$ for finite exponents $p$ and $q$, the analogous results were proved in \cite{TM5}. More specifically, it was proved that $V_g: \mathcal{F}_{\psi_m}^p \to \mathcal{F}_{\psi_m}^q$ for $0<p\leq q<\infty$ is   bounded if and only if  $g$ is a complex polynomial of degree not exceeding. Compactness was described in terms of the degree of $g$ being  at most  one. On the other hand,  if  $0<q<p<\infty$, then $V_g: \mathcal{F}_{\psi_m}^p \to \mathcal{F}_{\psi_m}^q$ is bounded(compact) if and only if  $g$ is a polynomial of degree not bigger than one. The same  conclusion as above holds with  respect to the spectrum of the operators $V_g$. Thus, our main results now could be seen as a completion of  the missing gap in \cite{TM5} when at least one of the exponent is at infinity. From this and the corresponding results in \cite{TM5}, we now, in addition, conclude that  the boundedness, compactness,  and spectrum of $V_g$ are independent of the order $m$ of the Fock--Sobolev spaces.

We may now  note  that  if we set $\psi_{(\alpha, t)}(z)= \alpha |z|^t, \ t>0, \alpha >0$, then the boundedness and compactness properties of   $V_g$ acting on the growth type weighted  spaces   $\mathcal{F}_{\psi_{(\alpha, t)}}^\infty$ has been described  in \cite{Bonet1} while its spectrum was identified later in \cite{Bonet2} \footnote{We would like to  thank the anonymous reviewer for bringing  the works in \cite{Bonet1} and \cite{Bonet2} to our attention}. In context of the  works in these two articles,  our results in Theorem~\ref{thm1} for the case when  $p= \infty$ can be considered as extension results obtained by  making logarithmic   perturbation of the weight function $\psi_{(\frac{1}{2}, 2)}(z)= \frac{1}{2}|z|^2$ into $\psi_m(z)= \frac{1}{2}|z|^2-m\log(1+ |z|) $ for all positive integers $m$.  It follows that the results are preserved under such perturbation.  We also  note in passing  that most of the techniques used  to prove our results in here are different from those used in \cite{Bonet1, Bonet2}.
\section{Preliminaries}
Note that for  each nonnegative integer $m$, the spaces $\mathcal{F}_{\psi_m}^2 $  are  reproducing kernel Hilbert
spaces with  kernel $K_{(w,m)}$ and normalized reproducing kernel functions $k_{_{(w,m)}}$  for a  point $w$ in $\CC.$  An explicit expression for
$K_{(w,m)}$ is still unknown. By  Proposition~2.7 of  \cite{CCK}, for each $w$ in $\CC$, we have the asymptotic relation
\begin{equation}
\label{asymptotic}
\|K_{_{(w,m)}}\|_{(2, m)}^2 \simeq  e^{2\psi_m(w)}.
\end{equation}
Observe that when $m=0,$ the space $\mathcal{F}_m^2$ reduces to the classical Fock  space $\mathcal{F}^2$, and in this  case  we  precisely have
 $\|K_{_{(z,0)}}\|_{(2, 0)}^2= e^{|z|^2}$ and $K_{(w,0)}(z)= e^{\overline{w}z}$. For other  values of $p$, by Corollary~14 of \cite{RCKZ},   we only have a one sided estimate
\begin{align}
\label{forall}
\|K_{(w, m)}\|_{(p,m)} \lesssim  e^{\psi_m(w)}.
\end{align}
  Studying  Volterra-type integral operators  in normed spaces  gets handy when  the norms in the target spaces of the operators are   described   in terms of Littlewood--Paley type  formula. These operators have been extensively studied in the spaces where such  formulas  are  happened to be known. Our next key lemma does this job by characterizing  the growth type  Fock--Sobolev spaces in terms of
 derivatives  which is interest of its own.
\begin{lemma}\label{lem1}
Let $f$ be an entire function on $\CC$. Then  $f$ belongs to the spaces $\mathcal{F}_{\psi_m}^\infty$ if and only if
\begin{align*}
\sup_{z\in \CC} \frac{|f'(z)|e^{-\psi_m(z)}}{1+\psi_m'(z)} <\infty,
\end{align*} and in this case we estimate the $\mathcal{F}_{\psi_m}^\infty$-norm of $f$ as
\begin{align}
\label{norm}
\|f\|_{(m, \infty)} \simeq |f(0)| + \sup_{z\in \CC} \frac{|f'(z)|e^{-\psi_m(z)}}{1+\psi_m'(z)}.
\end{align}
\end{lemma}
\emph{Proof}. For each function $f$ in $\mathcal{F}_{\psi_m}^\infty$, we make the estimate
\begin{align}
\label{part1}
|f(z)-f(0)| \leq \int_0^1 |z| |f'(tz)|dt\quad \quad \quad \quad \quad \quad \quad \quad \quad \quad \quad \quad \quad \quad \quad \quad \quad \quad \nonumber\\
= \int_{0}^1 \frac{|f'(tz)|e^{-\psi_m(t|z|)}}{1+\psi_m'(t|z|)} |z|(1+\psi'_m(t|z|))e^{\psi_m(t|z|)}dt  \quad \nonumber\\
\leq \Bigg(\sup_{z\in \CC} \frac{|f'(z)|e^{-\psi_m(z)}}{1+\psi_m'(z)}\Bigg) \int_{0}^1  |z|(1+\psi'_m(t|z|))e^{\psi_m(t|z|)}dt.
\end{align}
A straight forward integration by substitution shows that
\begin{align*}
\int_{0}^t |z|(1+\psi'_m(t|z|))e^{\psi_m(t|z|)}dt \lesssim e^{\psi_m(z)}.
\end{align*} Taking into account this  and  estimate \eqref{part1},   we obtain
\begin{align*}
\|f-f(0)\|_{(m, \infty)}  \lesssim \sup_{z\in \CC} \frac{|f'(z)|e^{-\psi_m(z)}}{1+\psi_m'(z)},
\end{align*} from which  and triangle inequality we deduce the one sided inequality
\begin{align}
\|f\|_{(m, \infty)}  \lesssim |f(0)| + \sup_{z\in \CC} \frac{|f'(z)|e^{-\psi_m(z)}}{1+\psi_m'(z)}.
\end{align}
 To prove the converse inequality, we first observe that by  subharmonicity of  $|f|$,
\begin{align}
\label{part2}
|f(z)| \lesssim \int_{D(z, 1)} | f(w)| dA(w) = \int_{D(z, 1)}  e^{\psi_m(w)} \big(|f(w)| e^{-\psi_m(w)}\big)dA(w),
\end{align} where $D(z, \delta)$ refers to a disk  of center at $z$ and radius $\delta$.
On the other hand,  for $w \in D(z, 1)$,  the estimate $\psi_m(w)\simeq \psi_m(z)$ holds. Thus, taking this into account in \eqref{part2} leads
\begin{align*}
|f'(z)| \lesssim \frac{d}{dz} \bigg(\int_{D(z, 1)}  e^{\psi_m(w)} \big(|f(w)| e^{-\psi_m(w)}\big)dA(w)\bigg)\quad \quad \quad \quad  \nonumber\\
\simeq \frac{d}{dz} e^{\psi_m(z)} \int_{D(z, 1)}  |f(w)| e^{-\psi_m(w)}dA(w)\nonumber\\
\lesssim e^{\psi_m(z)} \psi'_m(z) \|f\|_{(m, \infty)} \leq  \|f\|_{(m, \infty)}\frac{\big(1+\psi'_m(z)\big) }{e^{-\psi_m(z)}},
\end{align*} and from which we obtain  the remaining estimate
\begin{align*}
|f(0)| + \sup_{z\in \CC} \frac{|f'(z)|e^{-\psi_m(z)}}{1+\psi_m'(z)} \lesssim  \frac{|f(0)|e^{-\psi_m(0)}}{e^{-\psi_m(0)}} + \|f\|_{(m, \infty)} \lesssim \|f\|_{(m, \infty)}
\end{align*}  as desired and completes the proof.

As pointed,  the approximation formula \eqref{norm} is in the spirit of the famous  Littlewood--Paley  formula for entire functions in the growth type  space
$\mathcal{F}_{\psi_m}^\infty$. The corresponding  formula in  $\mathcal{F}_{\psi_m}^p$ for finite $p$ was obtained
in \cite{TM5} and reads as  \begin{align}
\label{Paley}
\|f\|_{(m,p)}^p \simeq |f(0)|^p + \int_{\CC} |f'(z)|^p \frac{e^{-p\psi_m(z)}}{(1+\psi_m'(z))^p} dm(z).
\end{align} Both formulas \eqref{norm} and \eqref{Paley} will be used repeatedly in our subsequent  considerations.
\begin{lemma}\label{lem2} Let $g(z)= az^2+bz+c$ and $|\lambda|>2|a|$. Then for any entire function $f$ on $\CC$ it holds that
\begin{align}
\label{smooth}
\sup_{z\in \CC} |e^{\frac{g(z)}{\lambda}} f(z)|e^{-\psi_m(z)} \lesssim  |f(0)| + \sup_{z\in \CC} \frac{|e^{\frac{g(z)}{\lambda}}||f'(z)|e^{-\psi_m(z)}}{1+\psi_m'(z)}.
\end{align}
\end{lemma}
This lemma provides a key tool to prove our main result on the spectrum of $V_g$ in the next section.\\
\emph{Proof}. Arguing as in the proof of the above lemma, we make the pointwise estimate
\begin{align*}
%\label{part11}
|f(z)-f(0)|| e^{\frac{g(z)}{\lambda}}| \leq e^{\frac{g(z)}{\lambda}}\int_0^1 |z| |f'(tz)|dt  \quad \quad \quad \quad \quad \quad \quad \quad \quad \quad \quad  \nonumber\\
\leq  \Bigg(\sup_{z\in \CC}\frac{ |e^{\frac{g(z)}{\lambda}}f'(z)|e^{-\psi_m(z)}}{1+\psi_m'(z)}\Bigg) \int_{0}^t |z|(1+\psi'_m(t|z|))e^{\psi_m(t|z|)}dt\nonumber\\
\lesssim  e^{\psi_m(z)}\Bigg(\sup_{z\in \CC}\frac{ |e^{\frac{g(z)}{\lambda}}f'(z)|e^{-\psi_m(z)}}{1+\psi_m'(z)}\Bigg).
\end{align*} From this and triangle inequality,  we deduce
\begin{align*}
\sup_{z\in \CC} |f(z) e^{\frac{g(z)}{\lambda}}|e^{-\psi_m(z)} \lesssim \sup_{z\in \CC} |f(0) e^{\frac{g(z)}{\lambda}}|e^{-\psi_m(z)}+ \sup_{z\in \CC} \frac{ |e^{\frac{g(z)}{\lambda}}f'(z)|e^{-\psi_m(z)}}{1+\psi_m'(z)}\nonumber\\
\leq |f(0)|+ \sup_{z\in \CC} \frac{ |e^{\frac{g(z)}{\lambda}}f'(z)|e^{-\psi_m(z)}}{1+\psi_m'(z)},
\end{align*} where for the last inequality we used the assumption that $|\lambda|>2|a|$ and
\begin{align}
\label{used}
\sup_{z\in \CC}|e^{\frac{g(z)}{\lambda}}|e^{-\psi_m(z)}= \sup_{z\in \CC} e^{\Re\big(\frac{az^2+bz+c}{\lambda}\big)-\psi_m(z)}\quad \quad \quad \quad \quad \quad \quad \quad \quad \nonumber\\
\lesssim \sup_{z\in \CC} e^{\big(\frac{|a|}{|\lambda|}-\frac{1}{2}\big)|z|^2+\frac{|bz|}{|\lambda|}+m\log(1+|z|)} \lesssim 1.
\end{align}
\section{Proof of the main result}
In this section we prove our main result.  We begin with the proof of part (a) of (i). If $g$ is a polynomial of degree not exceeding two, then $|g'(z)|\lesssim 1+ \psi'_m(z)$ for all $z\in \CC$.  Taking this into account and applying Lemma~\ref{lem1} leads to
\begin{align*}
\|V_gf\|_{(m, \infty)} \simeq \sup_{z\in \CC} \frac{|g'(z)||f(z)|e^{-\psi_m(z)}}{1+\psi_m'(z)} \leq \sup_{z\in \CC} \bigg(\frac{|g'(z)|}{1+\psi_m'(z)}\bigg)   \sup_{z\in \CC}\bigg( |f(z)|e^{-\psi_m(z)}\bigg)\nonumber\\
  \lesssim \|f\|_{(m, \infty)} \lesssim \|f\|_{(m, p)},\quad \quad \quad
\end{align*} where for the last inequality we used the monotonicity property $\mathcal{F}_{\psi_m}^p \subseteq \mathcal{F}_{\psi_m}^\infty$, and from which boundedness of $V_g$ follows.\\
 To prove the  converse,   for  each point  $w\in \CC^n$  we  consider   the  sequence of
 functions  $\xi_{(w, m)}(z)= e^{-\psi_m(w)} K_{(w,m)}(z).$ Then,
\begin{equation}\label{normb} \|\xi_{(w,m)}\|_{(m,p)} \lesssim 1 \end{equation}  independent of
$p$ and $w$ which follows from  \eqref{forall} for $p<\infty$ and  from a simple argument for $p= \infty.$
Applying $V_g$ to  such a sequence  yields
\begin{align*}
\|V_g\| \gtrsim \| V_g \xi_{(w, m)} \|_{(m,\infty)}\simeq
 \sup_{z\in \CC} \frac{|g'(z)||\xi_{(w, m)}(z)|e^{-\psi_m(z)}}{1+\psi_m'(z)}\nonumber\\
 \geq   \frac{|g'(z)||\xi_{(w, m)}(z)|e^{-\psi_m(z)}}{1+\psi_m'(z)}
\end{align*}
for all points $w$ and $z$ in $\CC.$ In particular,  setting  $w= z$  and applying \eqref{asymptotic}  gives
\begin{equation*}
\|V_g\| \gtrsim  \frac{|g'(w)||\xi_{(w, m)}(w)|e^{-\psi_m(w)}}{1+\psi_m'(w)}\simeq  \frac{|g'(w)|}{1+\psi_m'(w)},
\end{equation*}  which holds only  when  $|g'(w)|\lesssim 1+\psi_m'(w)$ for all $w$. This again  holds if and  only if $g$ is a complex polynomial of degree not exceeding two as asserted, and completes the proof of part (a).

To prove part (b) of (i), we  may first assume that $V_g$ is compact, and observe that the  sequence  $\xi_{(w, m)}$ converges to zero as
$|w| \to \infty,$ and the convergence is uniform on compact subset of $\CC.$ Then,  arguing as in the previous part, compactness of $V_g$ implies
\begin{align*}
0= \limsup_{|w|\to \infty} \| V_g \xi_{(w, m)} \|_{(m,\infty)}\simeq \limsup_{|w|\to \infty}\bigg( \sup_{z\in \CC} \frac{|g'(z)||\xi_{(w, m)}(z)|e^{-\psi_m(z)}}{1+\psi_m'(z)}\bigg)\nonumber\\
\gtrsim \limsup_{|w|\to \infty} \frac{|g'(w)||\xi_{(w, m)}(w)|e^{-\psi_m(w)}}{1+\psi_m'(w)}\simeq \limsup_{|w|\to \infty} \frac{|g'(w)|}{1+\psi_m'(w)},
\end{align*} from which it follows that $g$ is a polynomial of degree not exceeding one.\\
Conversely, assume that  $g(z)= az+b, \ a, b \in\CC$. Then,   obviously,  $g$ belongs to the space $\mathcal{F}_{\psi_m}^\infty$. We aim to show that $V_g: \mathcal{F}_{\psi_m}^p \to \mathcal{F}_{\psi_m}^\infty$ is compact.   To this end, let $f_j$ be a sequence of functions in $\mathcal{F}_{\psi_m}^p$ such that $\sup_j
\|f_j\|_{(m,p)}<\infty$ and $f_j$ converges uniformly to
zero on compact subsets of $\CC$ as $j\to \infty.$   Since $\psi'_m(z) \to \infty$ as $|z| \to \infty, $  for each $\epsilon >0$ there  exists
a positive $N_1$ such that
\begin{equation*}\frac{1}{1+\psi'_m(z)} <\epsilon
\end{equation*} for all $|z| > N_1.$  From this and  Lemma~\ref{lem1}, we obtain
 \begin{align*}
 |V_g f_j(z) | e^{-\psi_m(z)}\lesssim
 \frac{| a f_j(z)|e^{-\psi_m(z)}}{1+\psi'_m(z)}
  \lesssim  \|f_j\|_{(p,m)}\frac{ |a|} {1+\psi'_m(z)}
  \lesssim \epsilon
 \end{align*} for all $|z| > N_1$ and all $j.$ On the other hand,  if $|z|\leq N_1,$ then
\begin{align*}
 | V_g f_j(z) | e^{-\psi_m(z)}\lesssim
 \frac{| g'(z) f_j(z)|e^{-\psi_m(z)}}{1+\psi'_m(z)}
  \leq\|g\|_{(\infty,m)}\sup_{|z|\leq N_1 } |f_j(z)|\\
  \lesssim  \sup_{|z|\leq N_1 } |f_j(z)|  \to 0  \ \text{as}\ \ j\to \infty.
   \end{align*}
\emph{Part (ii)}: Applying the estimate in \eqref{Paley},  $V_g: \mathcal{F}_{\psi_m}^\infty \to \mathcal{F}_{\psi_m}^p$ is   bounded if and only if
the inequality
\begin{align}
\label{carleson}
\|V_gf\|_{(m, p)}^p= \int_{\CC} \frac{|f(z)|^p |g'(z)|^p }{(1+\psi'_m(z))^p} e^{-p\psi_m(z)} dA(z)\quad \quad \quad \quad \quad \quad  \quad \quad \quad \quad \quad\nonumber\\
= \int_{\CC} \frac{|f(z)|^p |g'(z)|^p (1+|z|)^{mp}}{(1+\psi'_m(z))^pe^{\frac{p}{2}|z|^2}}  dA(z)
 \lesssim \|f\|_{(m, \infty)}^p
 \end{align} holds for each $f\in  \mathcal{F}_{\psi_m}^\infty$. It means that if we set
\begin{align*}
d\mu_{(g, p)}(z)= \frac{|g'(z)|^p (1+|z|)^{mp}}{\big( 1+\psi_m'(z)\big)^p} dA(z),
\end{align*} then the inequality in \eqref{carleson} holds  if and only if $\mu_{(g, p)}$ is an $(\infty, p)$ Fock--Carleson measure.  Then an application of Theorem~2.4 of \cite{TM4} immediately gives that the statements in (a) and (b) are equivalent, and any of these holds if and only if $\tilde{\mu}_{(t, mp)}$  belongs to $L^1(\CC, dA)$  for some or any positive $t$ where
 \begin{align*}
\tilde{\mu}_{(t, mp)}(w)=\int_{\CC} \frac{e^{-\frac{t}{2}|z-w|^2}}{ (1+|z|)^{mp}} d\mu_{(g, p)}(z).
 \end{align*}
 Having singled out this equivalent reformulation, our next task will be to investigate further  the new formulation.
 Let us  first assume $\tilde{\mu}_{(p, mp)}$  belongs to $L^1(\CC, dA)$ and plan to show  $g$ is a complex polynomial of  degree  not exceeding one. Then,  using the fact that $1+\psi_m'(z) \simeq 1+\psi'_m(w)$ and  whenever $w$  belongs to the disk $D(z, 1)$,  and subharmonicity of $|g'|^p$, we infer
 \begin{align*}
 \frac{|g'(z)|^p}{(1+\psi'_m(z))^p} \lesssim \int_{D(z, 1)} \frac{|g'(w)|^p}{(1+\psi'_m(w))^p}  dA(w)\lesssim     \int_{D(z, 1)} \frac{e^{-\frac{p}{2}|z-w|^2}}{ (1+|w|)^{mp}}d\mu_{(g, p)}(w)\nonumber\\
  \lesssim  \tilde{\mu}_{(p, mp)}(z).
 \end{align*}
 Integrating the above shows that
 \begin{align*}
\int_{\CC}\frac{|g'(z)|^p}{(1+\psi'_m(z))^p}dA(z) \lesssim \int_{\CC}\tilde{\mu}_{(t, mp)}(w)dA(w) <\infty
 \end{align*} holds only if $g'$ is a constant for  $p>2$, and $g'= 0$ whenever $p\leq 2.$\\
 On the other hand if $g'= c= $constant, then
 \begin{align*}
\int_{\CC}\tilde{\mu}_{(p, mp)}(w) dA(w)=\int_{\CC}\bigg(\int_{\CC} \frac{e^{-\frac{p}{2}|z-w|^2}}{ (1+|z|)^{mp}} d\mu_{(g, p)}(z) \bigg)dA(w)\\
= \int_{\CC}  \int_{\CC} \frac{e^{-\frac{p}{2}|z-w|^2} |c|^p }{\big( 1+\psi_m'(z)\big)^p}  dA(z) dA(w)\nonumber\\
\simeq \int_{\CC} \frac{|c|^p}{\big( 1+\psi_m'(w)\big)^p} dA(w)<\infty,
 \end{align*} as $c= 0$ whenever $p\leq 2$ by our assumption,  and $\psi_m'$  is $L^p$ integrable for all $p>2$.  This  completes the proof of part (ii) of the main result.

 \emph{Part(iii)}:
We assume that $V_g$ is bounded on $\mathcal{F}_{\psi_m}^\infty$ and hence  $g(z)= az^2+ bz+c.$  Then,   by linearity of the integral we may write
$\lambda I- V_g= (\lambda I- V_{g_1})- V_{g_{2}}$ where $g_1(z)= az^2$ and $g_2(z)= bz+c$.  A simple analysis shows that $\lambda I- V_g$ and $ \lambda I- V_{g_1}$ are injective maps. On the other hand, by part (b) of (i) in the theorem, $V_{g_2}$ is a compact operator. Thus $\sigma(V_{g_2})= \{0\}.$ We shall proceed to consider the corresponding case with $g_1.$   We may first observe that  if  $\lambda\neq 0,$  then the equation  $\lambda f-V_g f= h$  has the unique analytic solution
\begin{align}
\label{first}
f(z)=(\lambda I- V_{g_1})^{-1}h(z)= \frac{1}{\lambda} h(0) e^{\frac{{g_1}(z)}{\lambda}} + \frac{1}{\lambda}  e^{\frac{{g_1}(z)}{\lambda}} \int_{0}^z  e^{-\frac{{g_1}(w)}{\lambda}} h'(w) dA(w).
\end{align}
 This can easily be  seen by solving an initial valued  first order linear ordinary differential equation
 $$ \lambda y'- {g_1}'y= h', \ \ \lambda f(0)= h(0). $$
 Recall that $(\lambda I- V_{g_1})^{-1}h=   R_{({g_1},\lambda)} h$
  is the Resolvent operator of $V_{g_1}$ at $\lambda$. It follows that $\lambda \in \CC$ belongs to the resolvent of $V_{g_1}$ whenever  $R_{({g_1},\lambda)}$ is a bounded operator. Since we assumed that $V_{g_1}$ is bounded and as  $\mathcal{F}_{\psi_m}^p $ contain the  constants, setting $h= 1$ in \eqref{first} shows that $R_{({g_1},\lambda)} 1= \frac{1}{\lambda}e^{g(z)/\lambda} \in \mathcal{F}_{\psi_m}^\infty $  for each $\lambda$ in the resolvent set of $V_{g_1}$.   On the other hand,  if $|\lambda| >2|a|$ then, from the estimation in \eqref{used} we have
   \begin{align*}
   \sup_{z\in \CC} |e^{{g_1}(z)/\lambda}|e^{-\psi_m(z)}<\infty,
   \end{align*} and from which we conclude that $|\lambda| >2|a|$ is a sufficient condition for the boundedness of $R_{({g_1},\lambda)}$ on $\mathcal{F}_{\psi_m}^\infty $. We aim to show that  the condition is  in fact  necessary as well.  To this end, let $f_1(z) =  \int_{0}^z e^{-\frac{g_1(w)}{\lambda}} f'(w)dA(w),$ and for   $2|a|<|\lambda|$,  Lemma~\ref{lem2} and then Lemma~\ref{lem1} implies
\begin{align*}
\|R_{(g_1,\lambda)}\|_{(m, \infty)} \leq \frac{|f(0)|}{|\lambda|}\| e^{\frac{g_1}{\lambda}} \|_{(m, \infty)} +  \frac{1}{|\lambda|}\sup_{z\in \CC} \frac{e^{\frac{g_1(z)}{\lambda}} |f'_1(z)|e^{-\psi_m(z)}}{1+\psi'_m(z)}\nonumber\\
\lesssim  \frac{|f(0)|}{|\lambda|} +  \frac{1}{|\lambda|}\sup_{z\in \CC} \frac{e^{\frac{g_1(z)}{\lambda}} |f'(z)|e^{-\psi_m(z)}}{1+\psi'_m(z)}
\lesssim  \frac{1}{|\lambda|}\| f \|_{(m, \infty)}.
\end{align*} We have now proved that for  a nonzero $\lambda$, the resolvent operator $R_{({g},\lambda)}$ is bounded if and only if $e^{{g}(z)/\lambda}$ belongs to $\mathcal{F}_{\psi_m}^\infty $, and this holds if and only if $|\lambda| >2|a|$.  From this our assertion
\begin{align*}
   \sigma(V_{g}) =\{0\}\cup \overline{ \{\lambda \in \CC\setminus\{0\}:e^{{g}(z)/\lambda} \notin \mathcal{F}_{\psi_m}^p\} }=\{\lambda \in \CC: |\lambda|\leq 2|a| \}
   \end{align*}immediately follows.


\begin{thebibliography}{BRSHZE}
\bibitem{Alman} A. Aleman, A class of integral operators on spaces of analytic functions,  Topics in
complex analysis and aperator theory,  3--30, Univ. M\'{a}laga,
M\'{a}laga, 2007.
\bibitem{Bonet1} J. Bonet and J. Taskinen,  Anote about Volterra operators on weighted Banach spaces of entire
functions, Math. Nachr. 288 (2015), 1216--1225.

\bibitem{Bonet2} J. Bonet, The spectrum of Volterra operators on weighted spaces of entire functions, Q. J. Math.
66 (2015), 799--807.

    \bibitem{9CH} S. Chandrasekhar, Radiative Transfer, Oxford Univ. Press, London, (1950).

\bibitem{CCK} H. R.  Cho, B. R. Choe, and  H. Koo, Linear combinations of composition operators on the Fock--Sobolev spaces, Potential Analysis, \textbf{4(41)} (2014), 1223--1246.
    \bibitem{RCKZ} R. Cho  and  K. Zhu,  Fock--Sobolev spaces and their Carleson measures, Journal of Functional Analysis
Volume,  263, Issue 8, \textbf{15}  (2012), 2483--2506.
\bibitem{Olivia2} O. Constantin and Ann-Maria  Persson, The spectrum of Volterra type integration operators on generalized  Fock spaces,
 Bull. London Math. Soc. (2015) doi: 10.1112/blms/bdv069.
\bibitem{Olivia} O. Constantin and Jos\'{e} \'{A}ngel Pel\'{a}ez, Integral Operators, Embedding Theorems and a Littlewood--Paley Formula on Weighted Fock Spaces,  Journal of Geometric Analysis, vol 26, 2(2016), 1109--1154
    \bibitem{Olivia1} O. Constantin, Volterra type integration operators on Fock spaces,  Proc. Amer. Math. Soc.,  \textbf{140}  (2012),
4247--4257.
\bibitem{3CC} C. Corduneanu, Integral Equations and Applications, Cambridge Univ. Press, Cambridge, (1991).

\bibitem{2KD} K. Deimling, Nonlinear Functional Analysis, Springer-Verlag, Berlin, (1985).
\bibitem{10HM} S. Hu, M. Khavanin and W. Zhuang, Integral equations arising in the kinetic theory of gases, Appl. Analysis,
\textbf{34} (1989), 201--266.

 \bibitem{TM3} T. Mengestie and S. I. Ueki,  Integral, Differential  and multiplication operators on weighted Fock spaces, Preprint, 2016.
 \bibitem{TM4} T. Mengestie, Carleson type measures for Fock--Sobolev spaces,
Complex Analysis and Operator Theory,  \textbf{8} (2014), no 6, 1225--1256.


\bibitem{TM5} T. Mengestie, Spectral properties  of Volterra-type integral operators on Fock--Sobolev  spaces, Preprint, 2016.

\bibitem{TM0} T. Mengestie, Volterra type and weighted composition operators on  weighted Fock spaces,  Integr. Equ. Oper.
Theory, \textbf{76} (2013),  81--94.
\bibitem{JPP} J. Pau and J. A. Pel\'{a}ez,  Embedding theorems and  integration operators on Bergman spaces  with rapidly decreasing weightes. J. Funct. Anal.,  \textbf{259 }(2010), 2727--2756.
    \bibitem{4DM} D. O'Regan and M. Meehan, Existence Theory for Nonlinear Integral and Integrodifferential Equations,
Kluwer Academic, Dordrecht, (1998).

  \bibitem{Si} A. Siskakis, Volterra operators on spaces of analytic functions-a survey,  Proceedings of the first advanced course
in operator theory and complex analysis, 51--68, Univ. Sevilla
Secr., Seville, 2006.

\end{thebibliography}
\end{document}